\documentclass[11pt, reqno, fleqn]{amsart} %\reqno = numero equations a droite

\usepackage{amsmath,amscd,amssymb,amsfonts,amsthm,courier,relsize,bm}
\usepackage{hyperref,enumitem,mathrsfs,slashed,mathtools, extarrows}
\setlist{leftmargin=10mm}
\usepackage[bottom]{footmisc}
\usepackage{xcolor}  %the next few lines set the colors of the hyperlinks, and remove the color boxes
\usepackage{geometry}
\usepackage[T1]{fontenc}
\usepackage{tikz}

\usepackage[bitstream-charter]{mathdesign} %In this one, replace \nu with \!\nu
\usepackage[T1]{fontenc}

\hypersetup{
    colorlinks,
    linkcolor={red!50!black},
    citecolor={blue!70!black},
    urlcolor={blue!80!black}
}

\newtheoremstyle{mytheoremstyle} % name
    {4mm}                    % Space above
    {4mm}                    % Space below
    {\itshape}                   % Body font
    {10mm}                           % Indent amount
    {\bfseries}                   % Theorem head font
    {. \, --- \,}                          % Punctuation after theorem head
    {0.5em}                       % Space after theorem head
    {}  % Theorem head spec (can be left empty, meaning ?normal?)

\theoremstyle{mytheoremstyle}
\newtheorem{theorem}{Theorem}[section]
\newtheorem*{theorem*}{Theorem}

\newtheorem{proposition}[theorem]{Proposition}

\newtheoremstyle{mydefinitionstyle} % name
    {4mm}                    % Space above
    {4mm}                    % Space below
    {}                   % Body font
    {10mm}                           % Indent amount
    {\bfseries}                   % Theorem head font
    {.\, --- \,}                          % Punctuation after theorem head
    {0.5em}                       % Space after theorem head
    {}  % Theorem head spec (can be left empty, meaning ?normal?)

\theoremstyle{mydefinitionstyle}    
\newtheorem{definition}[theorem]{Definition}

\newtheoremstyle{example}% name of the style to be used
  {4mm}% measure of space to leave above the theorem. E.g.: 3pt
  {4mm}% measure of space to leave below the theorem. E.g.: 3pt
  {}% name of font to use in the body of the theorem
  {10mm}% measure of space to indent
  {\bfseries}% name of head font
  {. \, --- \,}% punctuation between head and body
  { }% space after theorem head; " " = normal interword space
  {\thmname{#1}\thmnumber{ #2}\textnormal{\thmnote{ (#3)}}}

\theoremstyle{example}

\newtheorem*{example*}{Example}

\newtheorem*{examples*}{Examples}
\newtheorem{remark}[theorem]{Remark}

\newtheorem*{remark*}{Remark}
\newtheorem*{remarks*}{Remarks}

\usepackage{xpatch}
\xpatchcmd{\proof}{\hskip\labelsep}{\hskip6.5\labelsep}{}{}  %% change 5 here as you wish

\newcommand{\ignore}[1]{}

\newcommand{\ol}[1]{\overline{#1}}

\newcommand{\ti}[1]{\widetilde{#1}}
\newcommand{\wh}[1]{\widehat{#1}}

\renewcommand{\i}{{\mathrm{i}}}
\def\Ad{\ensuremath{\textnormal{Ad}}}

\def\g{\ensuremath{\mathfrak{g}}}

\def\E{\ensuremath{\mathscr{E}}}

\def\R{\ensuremath{\mathcal{R}}}
\def\S{\ensuremath{\mathfrak{S}}}

\def\L{\ensuremath{\mathscr{L}}} 
\def\f{\ensuremath{\mathfrak{f}}}

\def\Cl{\ensuremath{\textnormal{Cl}}}
\def\Cliff{\ensuremath{\textnormal{Cliff}}}

\def\C{\ensuremath{\mathbb{C}}}
\def\Z{\ensuremath{\mathbb{Z}}}
\def\K{\ensuremath{\mathbb{K}}}
\def\D{\ensuremath{\mathbf{D}}}

\def\i{\ensuremath{\mathrm{i}}}

\def\Hom{\ensuremath{\textnormal{Hom}}}

\def\index{\ensuremath{\textnormal{Index}}}
\def\supp{\ensuremath{\textnormal{supp}}}

\def\KK{\ensuremath{\textnormal{KK}}}

\def\K{\ensuremath{\textnormal{K}}}

\def\ch{\ensuremath{\textnormal{ch}}}

\usetikzlibrary{arrows,chains,matrix,positioning,scopes}
\makeatletter
\tikzset{join/.code=\tikzset{after node path={%
\ifx\tikzchainprevious\pgfutil@empty\else(\tikzchainprevious)%
edge[every join]#1(\tikzchaincurrent)\fi}}}
\makeatother
\tikzset{>=stealth',every on chain/.append style={join},
         every join/.style={->}}
\tikzstyle{labeled}=[execute at begin node=$\scriptstyle,
   execute at end node=$]
\begin{document}
\newgeometry{bottom=2.5cm, top=2cm, hmargin=2.5cm}
\sloppy
\title{A transverse Chern character for compact Lie group actions}

%\title{A Kasparov-like approach of the Berline--Vergne Chern character for transversally elliptic operators}

\author{Rudy Rodsphon}
\date{September 1, 2022}
\address{Washington University \\ 1 Brookings Dr \\ Department of Mathematics $\&$ Statistics, Cupples Hall \\ St Louis MO 63130, United States}
\email{rrudy@wustl.edu}

\maketitle

\vspace{-0.5cm}
\setlength{\parindent}{0mm}
\setlength{\mathindent}{15mm}

\begin{abstract}
This article aims to explore new perspectives offered by Kasparov's recent work on transverse index theory in the context of actions of compact Lie groups on a manifold.
\end{abstract}

\vspace{5mm}

\section*{Introduction} \setcounter{section}{0}

\vspace{2mm}
 
Transverse index theory for actions of compact Lie groups (on manifolds) consists in the study of the index problem for operators that are elliptic only in directions that are transverse to the orbits of the action. \\

A starting point is the monograph \cite{AtiTransEll}, which summarizes the investigations of Atiyah--Singer on the subject, and solves the problem in many interesting cases from a K-theoretic viewpoint, nonetheless without putting an end to the story, because of the lack of a general explicit formula. \ignore{In the eighties, Connes apparently communicated to Quillen that his superconnection formlism ought to be useful for the index problem.} In the nineties, Berline and Vergne \cite{BV1, BV2} solved the Atiyah--Singer transverse index problem, indeed through the use of superconnection techniques, which was refined and simplified by Paradan and Vergne in a series of papers \cite{PV1, PV2, PV3}. \\

The basic idea in these works is to start from the fact that an elliptic symbol may be viewed as superconnection form; as such, Quillen's superconnection methods provide a natural representative of its equivariant Chern character, which has a Gaussian shape peaking within the characteristic set of the symbol in question. When the operator is transversally elliptic, the rapid decay property does not hold anymore in the orbital directions. Berline--Vergne \cite{BV1} address this issue showing that for a certain class of transverse symbols, the use of the Liouville form and of the moment map on the cotangent bundle enables the addition of an oscillatory factor in the orbital directions, having the effect of restoring the rapid decay property of the Chern character, but only in the sense of distributions (over the Lie algebra). Paradan and Vergne \cite{PV2} make significant improvements, introducing another choice of additional factor that makes the Chern character have compact support, eliminating at the same time the need to restrict the construction to a subclass of transversally elliptic symbols. \\

Kasparov's recent work \cite{Kas16} takes a priori the subject in another direction via KK-theory, and obtains an index theorem which essentially states that the $K$-homology class of a transversally elliptic operator and its K-theory symbol class are dual in the sense of KK-theory. Heuristically, Kasparov's index formula is meant to reduce the calculation of a transverse index to a classical index, and consists essentially in adding the fundamental (orbital) vector field relative to the group action to compensate the defect of ellipticity of the transverse symbol. More precisely, this is done by introducing a new symbol class that takes the form of a Kasparov product between the initial transverse symbol and the fundamental vector field. \\

One may note many similarities between these two different approaches, especially in view of the objects involved. The main difference is that the Berline--Paradan--Vergne approach takes place at the cohomological level, whereas Kasparov works at the K-theoretic level. \\

This short note advocates that Kasparov's ideas may be utilized from a cohomological standpoint as well: the key point is to view a certain representative of the aforementioned new symbol class as a superconnection form, which in effect, almost brings us in the territory of classical index theory. We then verify on the example of the circle action on the plane, and of the zero-operator on the circle, that this procedure is consistent with the calculations of Atiyah--Singer and Berline--Paradan--Vergne. \\

Naturally comes the question of establishing a precise relationship between the methods of Kasparov and Berline--Paradan--Vergne; this will not be discussed here, but we leave this direction for a future work, or eventually an expanded version of this note. 

\vspace{2mm}
 
\subsection*{Plan of the article} Section 1 gives a short account of Kasparov's work \cite{Kas16} on transverse index theory, and provides a definition of the Chern character. Section 2 applies it to the example of the circle action on the complex plane.

\subsection*{Acknowledgement}  As usual in articles revolving around his work, I'm grateful to Gennadi Kasparov for having taught me so much over the years. The present note is mostly about reflections I've had through discussions with him.

\subsection*{Framework}
We consider throughout the paper a complete Riemannian manifold $(M, g)$ (without boundary) equipped with an (isometric) action of a compact Lie group $G$ (the results of Section 1 also hold if $G$ is  locally compact and acts properly, isometrically). Let $\tau=TM$ and $\g=\mathrm{Lie}(G)$ denote respectively the tangent bundle of $M$ and the Lie algebra of $G$. We will not make much distinctions between $TM$ and the cotangent bundle $T^*M$, and will most of the time identify them tacitly via the metric. Similarly, $\g$ and its dual $\g^*$ are most of the time identified via a $G$-invariant product on $\g$. If $A, B$ are $C^*$-algebras, $\KK_{\bullet}^G(A,B)$, $K_{\bullet}^G(A)$ and $K^{\bullet}_G(B)$ the associated $G$-equivariant KK-theory, K-theory, K-homology groups. If moreover $A, B$ are $C_0(X)$-algberas,  $\R\KK_{\bullet}^G(X; A,B)$ denotes the representable KK-theory group of $A$ and $B$. 

%\vspace{2mm}

\section{Kasparov's approach of transverse index theory} \setcounter{section}{1} \label{Section1}

In this section, we review Kasparov's KK-theorerical approach of transverse index theory for Lie group actions \cite{Kas16} and the associated index theorem. 

\subsection{Orbital Clifford algebra}

Let $(M, g)$  be a complete Riemannian manifold (without boundary) equipped with an (isometric) action of a compact Lie group $G$. For every $x \in M$, the derivative of the action of $G$ on $x$ reads:   
\[ \rho_x: v \in \g \longmapsto \left. \dfrac{d}{dt}\right\vert_{t=0} \exp(-tv) \cdot x \, \in T_x M \]
and yields a smooth bundle map $\rho: \g_M:=M \times \g \to TM$.  The \emph{orbital tangent field} $\Gamma \subset TM$ is the continuous field of subspaces defined as the image $\Gamma = \mathrm{Im}(\rho)$ of $\rho$. The family of Clifford algebras 
\[ \Cliff(\Gamma) = \bigsqcup_{x \in M} \Cliff(\Gamma_x) \subset \Cliff(TM) \] 
inherits in turn a continuous field structure from $\Gamma$. The \emph{orbital Clifford algebra} $\Cl_\Gamma(M)$ is the $C^\ast$-algebra generated by the $C_0$-sections of the continuous field $\Cliff(\Gamma)$ over $M$ (equipped with the sup-norm). Up to isomorphism, this definition if independent of the choice of the metric on $M$. \\

On the other hand, there is a natural extension of the $G$-action on $M$ to $\g_M$ given by $g \cdot (x,v) = ( g \cdot x, \Ad(g)v )$, for which $\rho$ is $G$-equivariant. This action being proper, one can equip $\g_M$ with a $G$-invariant Riemannian metric on the bundle $\g_M$. We can then consider the transpose operator $\rho^t \colon TM \rightarrow \g_M$, and introduce a \emph{smooth} bundle map $\varphi \colon TM \rightarrow TM$ to be the composition $\varphi=\rho \rho^t$, which is roughly speaking the projection onto $\Gamma$: indeed, note that \emph{a vector $\xi \in T_xM$ is orthogonal to $\Gamma_x$ if and only if $\varphi_x(\xi)=0$}. 
\vspace{2mm}

\subsection{Transverse index and Dirac elements} Let $A$ be a $G$-equivariant pseudodifferential operator with symbol $\sigma_A$ acting on sections of a $G$-equivariant Hermitian vector bundle $E$. Recall that the operator $A$ is said to be \emph{transversally elliptic} if $\supp(\sigma_A)\cap T_G^*M$ is compact, where $T_G^\ast M=\{(x,\xi) \in T^*M \, ; \, \langle \xi, \Gamma_x \rangle = 0 \}$, and $\supp(\sigma_A)$ is the support of the symbol $\sigma_A$ (i.e the subset of $T^\ast M$ where $\sigma_A$ fails to be invertible). Whence it defines naturally a K-theory class $[\sigma_A]  \in \K^0_G(T_G^*M)=\K^G_0(C_0(T_G^*M))$. \\

When $M$ is compact, the restriction $A_\pi$ ($\pi \in \wh{G}$) of $A$ to each isotypical component is Fredholm \cite{AtiTransEll}, hence $A$ has a well-defined equivariant index 
\[ \index_G(A)=\sum_{\pi \in \wh{G}} \index(A_\pi)\pi \in \widehat{R}(G)=\Z^{\wh{G}}.\]
taking values in the generalized character ring $\widehat{R}(G)$. Form a K-homological standpoint, this fact means the following:

%When $M$ is not compact, Atiyah defines the index of $A$ by a reduction to the compact case, via an embedding of a relatively compact neighborhood $U$ of $\supp(\sigma_A) \cap T^*_GM$ into a compact manifold. We will take an alternate route, via KK-theoretical methods. \\

\begin{proposition} \emph{(cf. for example \cite[Proposition 6.4]{Kas16})} \, The pair $(L^2(M,E), A)$ induces a $K$-homology class 
\[[A] \in K^0(G \ltimes C_0(M)) \]
that we call the \emph{transverse index of $A$}. 
\end{proposition}

If $M$ is compact, then one can also view $[A]$ as a class in the K-homology group $K^0(C^*(G))$ by crushing $M$ to a point, and the Peter-Weyl theorem shows that $[A]$ coincides with Atiyah's index. \\ 

Suppose now that $E$ is a $G$-equivariant ($\Z_2$-graded) Clifford module bundle over $M$, and let $D$ be an associated Dirac-type operator. As it is elliptic, it induces a canonical $K$-homology class $[D] \in K^0(C_0(M))$ represented by the $K$-cycle $\left(L^2(M,E), F := D (1+D^2)^{-1/2} \right)$. The same K-cycle can be used to describe another $K$-homology class that takes into account the transverse geometry given by the $G$-action. 

\begin{definition} \label{thm:FundClass} (cf. for example \cite{Kas16, LRS19})

\begin{enumerate}[leftmargin=5mm]
\item[(i)]  The pair $(L^2(M,E), F)$ determines a $K$-homology class 
\[ [D_{M,\Gamma}] \in K^0(G \ltimes \Cl_\Gamma(M)), \] where the crossed product $G \ltimes \Cl_\Gamma(M)$ acts on $L^2(M,E)$ by multiplication (and convolution in $G$). This class will be referred to as the \emph{transverse Dirac element}. \\

\item[(ii)] When $E$ is the (complexified) exterior algebra bundle $\Lambda_{\C}^{\bullet} TM := \Lambda^{\bullet} TM \otimes \C$, with $D$ being the canonical Dirac-type operator associated to the de Rham differential $d$, the class $[D_{M,\Gamma}]$ promotes to a $K$-homology class 
\[[d_{M,\Gamma}] \in  K^0(G \ltimes \Cl_{\tau \oplus \Gamma}(M)), \] 
where $\Cl_{\tau \oplus \Gamma}(M) := \Cl_\tau(M) \otimes_{C_0(M)} \Cl_\Gamma(M)$. This class will be referred to as the \emph{transverse de Rham-Dirac element}.
\end{enumerate}
\end{definition}

The action of $\Cl_{\tau \oplus \Gamma}(M)$ on $L^2(M, \Lambda_{\C}^{\bullet} TM)$ is given, on real (co)vectors, by
\[ \xi_1 \oplus \xi_2 \longmapsto \mathrm{ext}(\xi_1) + \mathrm{int}(\xi_1) + \i (\mathrm{ext}(\xi_2) - \mathrm{int}(\xi_2)) \]
where $\xi_1, \xi_2$ are respectively sections of $TM$ and $\Gamma$. \\

The fact that one indeed gets K-homology classes this way is not completely obvious, especially the verification that the space $[D_{M, \Gamma}, G \ltimes \Cl_\Gamma(M)]$ consists of compact operators; see \cite{Kas16} for (ii), and \cite{LRS19} for (i). %The general idea is that the $G$-action by convolution allows to make up for the unboundedness of the elements in $[D, \Cl_\Gamma(M)]$.  

%Another alternative way to proceed is to observe that $[D_{M, \Gamma}] \in K^0(G \ltimes \Cl_\Gamma(M))$ is the image of $[D] \in K^0(G \ltimes C_0(M))$ under the projection map  $\Cl_{\Gamma}(M) \to C_0(M)$. 

\vspace{2mm}

\subsection{Transversally elliptic symbols}\label{sec:symbalg}

The precise relationship between the classes $[\sigma_A] \in\K^0_G(T_G^*M)$ and $[A] \in K^0(G \ltimes C_0(M))$ is the object of a highly non-trivial index theorem of Kasparov that we describe in next subsection. Before that, a replacement of the algebra $C_0(T_G^*M)$ is needed. 

\begin{definition} \, (\cite{Kas16}, Definition-Lemma 6.2) \, \label{def:symbalg}
The \emph{symbol algebra} $\S_\Gamma(M)$ is the norm-closure in $C_b(T^*M)$ (the algebra of continuous bounded functions on $T^*M$) of the set of all smooth, bounded functions $b(x,\xi)$ on $T^*M$, which are compactly supported in the $x$-variable, and satisfy (a) together with either (b) or (c), which are equivalent under (a):
\begin{enumerate}
\item The exterior derivative $d_xb(x,\xi)$ in $x$ is norm-bounded uniformly in $\xi$, and there is an estimate $|d_\xi b(x,\xi)|\le C(1+ \Vert \xi \Vert)^{-1}$ for a constant $C$ which depends only on $b$ and not on $(x,\xi)$.

\item The restriction of $b$ to $T_G^*M$ belongs to $C_0(T_G^*M)$. 

\item For any $\varepsilon >0$ there exists a constant $c_\varepsilon>0$ such that
\[ |b(x,\xi)| \leq c_\varepsilon \frac{1+ \Vert \varphi_x(\xi)\Vert^2}{1+ \Vert \xi \Vert^2}+\varepsilon, \qquad \forall x \in M, \xi \in T_xM.\]
\end{enumerate}
\end{definition}

Loosely speaking, item (c) says that $\frak{S}_\Gamma$ consists of  symbols with negative order in the transverse directions. \\

Let $\pi: T^*M \to M$ denote the canonical projection. Given a $G$-equivariant $\Z_2$-graded Hermitian vector bundle $E$, we can similarly define a Hilbert $\S_\Gamma(M)$-module, denoted $\S_\Gamma(E)$, as the norm-closure in the space of bounded sections of the pull-back bundle $\pi^\ast E$ satisfying similar conditions to those in Definition \ref{def:symbalg} (using the norm on the fibres of $p^\ast E$ induced by the Hermitian structure). \\

From now on, we refer to transversally elliptic operators (or symbols) according to the following definition.

\begin{definition} Let $A$ be a properly supported, odd, self-adjoint $G$-invariant pseudodifferential operator of order 0 acting on sections of a $G$-equivariant $\Z_2$-graded Hermitian vector bundle $E$. We will say that $A$ (or its symbol $\sigma_A$) is transversally elliptic if for every $a \in C_0(M)$,  $a\cdot(1-\sigma_A^2)\in \S_\Gamma(M)$. 
\end{definition}

Therefore, a transversally elliptic symbol naturally determines a class
\[ [\sigma_A]=[(\S_\Gamma(E),\sigma_A)]\in \R\KK^G(M;C_0(M),\S_\Gamma(M)).\]
By construction there is a $\ast$-homomorphism $\S_\Gamma(M)\rightarrow C_0(T_G^*M)$, hence a map
\[ \R\KK^G(M;C_0(M),\S_\Gamma(M))\rightarrow \R\KK^G(M;C_0(M),C_0(T_G^*M)).\]
In this sense the element $[\sigma_A]\in \R\KK^G(M;C_0(M),\S_\Gamma(M))$ can be viewed as a `refinement' of the `naive' class in $\R\KK^G(M;C_0(M),C_0(T_G^*M))$ defined by the symbol. 

\vspace{2mm}

\subsection{Kasparov's index theorem} To state Kasparov's index theorem relating the classes  $[A] \in K^0(G \ltimes C_0(M))$ and $ [\sigma_A]\in \R\KK^G(M;C_0(M),\S_\Gamma(M))$, we introduce the $C^*$-algebra
\[ \Cl_{\Gamma}(TM) := C_0(TM) \otimes_{C_0(M)} \Cl_{\Gamma}(M). \]
which is KK-equivalent to the symbol algebra $\S(M)$ through the KK-class described in the definition below. 

\begin{definition} \label{def:transell} %\cite[pp.1344--1345]{Kas16}
The element $[\f_{M,\Gamma}]\in \R\KK^G(M;\S_\Gamma(M),\Cl_\Gamma(TM))$ is the class represented by the pair $(\Cl_\Gamma(TM),\f_{M,\Gamma})$ where at a point $(x,\xi) \in T_xM$, the operator $\f_{M,\Gamma}(x,\xi)$ is left Clifford multiplication by $-\i \varphi_x(\xi)(1+ \Vert \varphi_x(\xi) \Vert^2)^{-1/2}$.
\end{definition}

This definition is in fact one of the main reasons to introduce the symbol algebra $\frak{S}_{\Gamma}(M)$. Being the symbol class of an orbital Dirac element, the class $[\f_\Gamma]$ should be thought of as an orbital Bott element. \\[-2mm]

Furthermore, notice that the $C^*$-algebras $\Cl_{\Gamma}(TM)$ and $\Cl_{\tau \oplus \Gamma}(M)$ are KK-equivalent. Indeed, recall that the $C^*$-algebras $C_0(TM)$ are KK-equivalent via an element $[d_\xi]\in \R\KK^G(M;C_0(TM),\Cl_\tau(M))$ referred to as the \emph{fiberwise Dirac element}. It is represented by a family of Dirac operators acting on the fibres $TM$, i.e the pair $\left( (L^2(T_xM) \otimes \Cl_{\tau_x})_{x \in M}, (F_x = D_x (1+D^2_x)^{-1/2})_{x \in M} \right)$, where for each $x \in M$, $D_x$ is a Dirac operator acting on the fiber $T_xM$. Its inverse can be described explicitely via a familly of fiberwise Bott elements. The element 
\[ [d_\xi] \otimes_{C_0(M)} 1_{\Cl_\Gamma(M)} \in \R\KK^G(M;\Cl_\Gamma(TM),\Cl_{\tau \oplus \Gamma}(M)) \]
then implements a KK-equivalence between $\Cl_{\Gamma}(TM)$ and $\Cl_{\tau \oplus \Gamma}(M)$. Motivated by the fact that the Dirac element on $TM$ induced by the Dolbeault operator splits as a product between $[d_\xi]$ and the Dirac element induced by the de Rham operator, Kasparov makes the following definiton. 

\begin{definition} \textbf{} We consider the following $\KK$-classes:
\begin{enumerate}[leftmargin=7mm] 
\item[(i)] The \emph{fiberwise transverse Dirac element} is the class $[d_{\xi, \Gamma}]$ obtained via the KK-product:
\[ [d_{\xi, \Gamma}] := [\f_{\Gamma}] \otimes_{\Cl_{\Gamma}(TM)}([d_\xi] \otimes_{C_0(M)} 1_{\Cl_\Gamma(M)})   \in \R\KK^G(M; \S_{\Gamma}(M),\Cl_{\tau \oplus \Gamma}(M))) \]

%\item[(ii)] The \emph{transverse Dolbeault element} is the class $[d_{\xi, \Gamma}]$ obtained via the KK-product:
%\[ [\overline{\partial}_{M, \Gamma}] := j^G ([d_{\xi, \Gamma}]) \otimes_{G \ltimes Cl_{\tau \oplus \Gamma}(M)} [d_{M, \Gamma}]  \in K_G^0(G \ltimes \S_{\Gamma}(M)) \]
%where $j^G$ is the descent map. \\[-2mm]

\item[(ii)] The \emph{transverse Dolbeault element} is the class $[\overline{\partial}_{TM,\Gamma}]$ obtained as the KK-product
\[ [\overline{\partial}_{TM,\Gamma}]=j^G([d_\xi]\wh{\otimes}_{M}1_{\Cl_\Gamma(M)})\wh{\otimes}_{G\ltimes \Cl_{\tau\oplus\Gamma}(M)}[d_{M,\Gamma}] \in \KK(G\ltimes \Cl_\Gamma(TM),\C).
\]%\footnote{Recall that for every $G$-$C^*$-algebras $A$ and $B$, the descent map $j^G: KK^G(A,B) \to KK(G \ltimes A, G \ltimes B)$.}.
\end{enumerate}
\end{definition}
Another way to represent the class $[\ol{\partial}_{TM, \Gamma}]$ is to use the Dolbeault operator $\ol{\partial}$ on the almost complex manifold $TM$. The latter induces a K-homology class via the operator $\ol{\partial}+\ol{\partial}^*$ via Definition \ref{thm:FundClass}(i). Proving that this class factors into the product above is overall similar to the elliptic case without group action, because the operators involved in the cycles are the same (cf. \cite{Kas88} for example). The only part that requires some attention, when checking is an $[d_{M,\Gamma}]$-connection, is the compactness of the commutators $[\ol{\partial}+\ol{\partial}^*, \Cl_{\tau \oplus \Gamma}]$, but this is done in \cite[Lemma 8.8]{Kas16}. \\ %the positivity condition of the Connes--Skandalis criterion,   %that the class of $\ol{\partial}+\ol{\partial}^*$ 

Before stating Kasparov's index theorem, it will be convenient to define another version of the symbol class.

\begin{definition} Let $A$ be a transversally elliptic operator, and let $[\sigma_A] \in \R\KK(M; C_0(M), \S_{\Gamma}(M))$ be its standard symbol class. The \emph{tangent Clifford symbol class} of $A$ is the element $[\sigma_A^{\mathrm{tcl}}]$ obtained as the KK-product:
\[ [\sigma_A^{\mathrm{tcl}}] := [\sigma_A] \otimes_{\S_{\Gamma}(M)} [\f_{\Gamma}] \in \R\KK(M, C_0(M), \Cl_{\Gamma}(TM)) \]
\end{definition}
Heuristically, the role of this product is to complete $[\sigma_A]$ into an elliptic symbol, so its index may be recovered from a pairing with the Dolbeault element. We can now state the following result of Kasparov, which converts this vague idea into a theorem. 

\begin{theorem} %\emph{(Kasparov, \cite[Theorem 8.18]{Kas16})} 
\label{thm:KaspInd}
Let $(M, g)$ be a complete Riemannian manifold (without boundary) equipped with a proper and isometric action of a Lie group $G$. Let $A$ be a transversally elliptic operator of order 0 on $M$, with symbol $\sigma_A$. Then
\[ [A]=j^G[\sigma_A^{\mathrm{tcl}}] \otimes_{G\ltimes \Cl_\Gamma(TM)}[\overline{\partial}_{TM,\Gamma}] \in \K^0(G\ltimes C_0(M)),\]
where $j^G$ denotes the descent map. 
\end{theorem}
%\[ [A] = j^G [\sigma_A] \otimes_{\S(M)} [\overline{\partial}_{M, \Gamma}] \in K^0(G \ltimes C_0(M)) \]
When $\sigma_A$ is transversally elliptic in Atiyah's sense, note that $[\sigma_A^{\mathrm{tcl}}]$ can also be seen as a class in $K_0(\Cl_{\Gamma}(TM))$. 

\begin{remark} \label{rk: trans to classical} To end this paragraph, note that the formula above is consistent with the classical index theorem when $[\sigma_A] \in \K_0^G(C_0(TM)))$ is the symbol class of an elliptic operator. Indeed, consider the canonical injection $\iota: C_0(TM) \to \Cl_{\Gamma}(TM)$. This induces a map $\iota_*: \K_0^{G}(C_0(TM)) \to \K_0^G(\Cl_{\Gamma}(TM))$. By associativity of the Kasparov product, we then have:
\[ j^G(\iota_*[\sigma_A]) \otimes [\overline{\partial}_{TM,\Gamma}] = j^G[\sigma_A] \otimes j^G[\iota] \otimes [\overline{\partial}_{TM,\Gamma}] =  j^G[\sigma_A] \otimes [\overline{\partial}_{TM}], \] 
where $[\iota] \in \KK(C_0(TM),\Cl_{\Gamma}(TM))$ is the KK-class associated to the injection $\iota$, and the K-homology class $[\overline{\partial}_{TM}] \in \K^0(G \ltimes C_0(TM))$ is the fundamental class induced by the Dolbeault operator on $TM$.  
\end{remark}
\vspace{1mm}

\section{Kasparov's fixed point formula}

This section is mostly expository and explains how Kasparov uses his machinery to tackle concrete calculations; it mostly follows  the work of Atiyah \cite{AtiTransEll}, while also bringing a handful of conceptual simplifications. The material is based on discussions and communications I've had with Kasparov; it may appear with more details and a better polished version in future works of his. 

%The goal of this section is to translate Kasparov's calculations into a cohomological approach on an example. We will keep the same notations as in previous section. 

\subsection{Reduction to toral actions}

%We start by summarizing how Kasparov explains the reduction of its transverse index theorem to the classical index theorem. Part of it also relies on observations of Atiyah, so we begin by translate these into Kasparov's framework. \\

The first step consists in a reduction to toral actions \cite[Chapter 4]{AtiTransEll}. Assume that $G$ is connected, let $T$ be the maximal torus of $G$ and let $\mathfrak{t}=\text{Lie}(T)$. With a slight abuse of notation, we still use $\Gamma$ to denote the orbital field in $M$ generated by the action of $T$.

\begin{proposition}
$K_0^{G} (\Cl_{\Gamma}(TM))$ is isomorphic to $K_0^{T}(\Cl_{\Gamma}(TM))$.  
\end{proposition}

Indeed, consider Kasparov's induction morphism $\mathrm{Ind}^{G}_{T}$, which gives a map: 
\[ K_0^{T}(\Cl_{\Gamma}(TM)) \longrightarrow KK_0^{G} \big(C(G/T),  \Cl_{\Gamma}(T(G \times_T M )) \big) \]

Then, observing that $G/T = T_G(G/T)$, and keeping in mind that the zero operator on $G/T$ is $G$-transversally elliptic, we compose the previous map with the (right) KK-product by $[\bm{0}] \in K^0(T_G^{*}(G/T))$, which yields a homomorphism:
\[ K_0^{T}(\Cl_{\Gamma}(TV)) \longrightarrow K_0^{G}\big(\Cl_{\Gamma}(T(G\times_T M))\big) \simeq K_0^{G}\big(\Cl_{\Gamma}(T(G/T \times M))\big) \]
where the last equality uses that $G\times_T M$ and $G/T \times M$ are $G$-equivariantly diffeomorphic via the map $(g,x) \mapsto (gT, g \cdot x)$. Finally, since $G/T$ is a rational variety, it comes with a Dolbeault element that has index $1$ (by a theorem of  Dolbeault). So $K_0^{G}(\Cl_{\Gamma}(T(G/T \times M)) \simeq K_0^{G}(\Cl_{\Gamma}(TM))$, and we get a map: %$[\ol{\partial}_{G/T}] \in K^{0}_G(C(G/T))$\\
\begin{align*}
K_0^{T}(\Cl_{\Gamma}(TM)) \longrightarrow K_0^{G} (\Cl_{\Gamma}(TM))
\end{align*}
which is shown to be an isomorphism by Atiyah \cite[Section 4]{AtiTransEll}. Note that the significance of $\Cl_{\Gamma}(TM)$ differs on both sides: on the right, the field $\Gamma$ is relative to the action of $T$, whereas the left considers the $G$-action. With this in mind, Atiyah's argument can be readapted straightforwardly in this context, and following his lead, we still denote the homomorphism just constructed by
\[ \mathrm{Ind}^{G}_{T}: K_0^{T}(\Cl_{\Gamma}(TM)) \longrightarrow K_0^{G} (\Cl_{\Gamma}(TM)). \]

The next point is a variant of the Atiyah--Segal localization method: one localizes the $R(T)$-module $K_0^{T}(\Cl_{\Gamma}(TM))$ with respect to the multiplicative family $S$ generated by elements $(1-e^{\i \theta}) \in R(T)$, where $\theta \mapsto e^{\i \theta}$ is any closed 1-parameter subgroup of $T$. As any proper subgroup $H \subset T$ is annihilated by a character of the latter type, the localization procedure only retains the contribution of the fixed point set $M^T \subset M$, so the restriction map (induced by the latter inclusion)
\[ K_0^T(\Cl_{\Gamma}(TM)) \longrightarrow K_0^T(\Cl_{\Gamma}(TM^T))\]
becomes an isomorphism after applying localization. \\

The normal bundle $N$ of the inclusion $M^T \hookrightarrow M$ naturally admits a complex structure because of the $T$-action. Identifying $N$ with a tubular neighborhood $U$ of $M^T$ within $M$, excision provides an isomorphism between the localized K-groups:
\[ K_0^T(\Cl_{\Gamma}(TU))_{\mathrm{loc}} \longrightarrow K_0^T(\Cl_{\Gamma}(TM^T))_{\mathrm{loc}}. \]
so we can replace $M$ by $U$ whenever it is convenient to do so. 

\subsection{Reduction to the classical index theorem}

Following up the preparations made in the previous subsection, we now replace the bundle $E \to M$ with $\mathcal{E} := E \times \Lambda^{\bullet} N \to M$ and consider its pull-back $p^*\mathcal{E} \to TM$. On the other hand, we replace the symbol $\sigma_A$ with the symbol $\ti{\sigma}_A := \sigma_A \otimes 1$. \\

Working out Kasparov's theorem \ref{thm:KaspInd} in this case, we need to deal with the symbol class
\[ [\ti{\sigma}_A^{\mathrm{tcl}}] := [\sigma_A] \otimes 1   \otimes_{\S_{\Gamma}(M)} [\f_{\Gamma}] \in \R\KK(M, C_0(M), \Cl_{\Gamma}(TM)) \]
represented by the pair 
\[ \big(C_0(TM, \E) \otimes_{C_0(TM)} \Cl_{\Gamma}(TM), \sigma_A \otimes 1 \otimes 1 + 1 \otimes 1 \otimes c(\varphi_x(\xi)) \big). \]
Performing a rotation homotopy in the two last tensors factors of the operator, we see (at least formally) that the symbol class $[\ti{\sigma}_A^{\mathrm{tcl}}]$ is also represented by the pair
\[ \big(C_0(TM, \E) \otimes \Cl_{\Gamma}(TM), \big(\sigma_A \otimes 1 + 1 \otimes c(\varphi_x(\xi))\big)\otimes 1 \big). \]
(up to weights involved in the KK-product, but this will not be needed). 

Now, the main point is that the symbol $\big(\sigma_A \otimes 1 + 1 \otimes c(\varphi_x(\xi))\big)$ is now elliptic, and by Remark \ref{rk: trans to classical}, so one may apply the usual index theorem to it. Applying the techniques of Atiyah--Segal, Kasparov gets the following result:

\begin{theorem}
The localized index of $A$ may be calculated via the localized $K$-theory class:
\[ \dfrac{[\left.\sigma_A^{\mathrm{tcl}}\right|_{TM^T}]}{\sum_{k\geq 0} (-1)^k [\Lambda^k p^*N]} \in \K^{T}_{0}(\Cl_{\Gamma}(TM^T))_{\mathrm{loc}} \]
where $[\Lambda^k p^{*}N]$ is the k-th exterior power of the pull-back bundle $p^*N$ over $TM^T$.   
\end{theorem} 

\begin{remark} If $N$ is a spin vector bundle, then one can also replace $\Lambda^{\bullet} p^*N$ with the spinor bundle of $N$. 
\end{remark}  
  
\section{A Chern character for transversally elliptic symbols and an example}

The goal of this section is to translate Kasparov's calculations into a cohomological approach on an example. We will keep the same notations as in previous section. 

\subsection{General discussion}
 
The basic idea is to view the operator $\big(\sigma_A \otimes 1 + 1 \otimes c(\varphi_x(\xi))\big)$, which represents the transversally elliptic symbol class $[\sigma_A^{\mathrm{tcl}}]$, as a superconnection form, after having substituted $E$ by its tensor product with a vector bundle $W$ that may be the exterior power/spinor bundle of the normal bundle $N$, or any other bundle in which $\Cliff \, \Gamma$ embeds into, in a way that is consistent with the action of $G$. A prior reduction to a toral action is not required. \\ %Notice that the latter is an elliptic symbol on $\E$, so from Remark \ref{rk: trans to classical}, the transverse index associated to the symbol class $[\ti{\sigma}_A^{\mathrm{tcl}}] \in K_0^T(\Cl_\Gamma(TM))$ may be calculated as the classical index associated to the symbol class $[\sigma_A \otimes 1 + 1 \otimes c(\varphi_x(\xi))] \in K_0^T(C_0(TM))$. \\

First, let us give some reminders about the Cartan model for the $G$-equivariant cohomology of $M$. In this setting, the latter is calculated through the cohomology of the complex $\big((\C[\g] \otimes \Omega^{\bullet}(M))^G, d_{\g} \big)$. The space $(\C[\g] \otimes \Omega^{\bullet}(M))^G$ is the subalgebra of $G$-invariant elements in the space of $\C[\g]$-valued differential forms, and the differential $d_{\g}$ is defined as 
\[ (d_{\g} \alpha)(v) = (d - \iota_v)(\alpha(v)) \quad ; \quad \forall v \in \g, \, \forall \alpha \in \C[\g] \otimes \Omega^{\bullet}(M).\] 
Here, $\iota_v$ denotes the contraction with the fundamental vector field $\textstyle x \in M \mapsto \rho_x(v)$. The Cartan homotopy formula for the Lie derivative $\L_v$ implies that $(d^2_{\g}\alpha)(v) = -\L_v \alpha(v)$, showing that $d_{\g}$ is a differential when restricting $\C[\g] \otimes \Omega^{\bullet}(M)$ to the subalgebra of $G$-invariant elements. \\  

Along the same lines, if $\D$ is a superconnection on a $\Z_2$-graded vector bundle $\E$ (acting on $\Omega(M,\E)$), then one defines its equivariant version acting on $\C[\g] \otimes \Omega(M, \E)$:
\[ \mathbf{D}_{\g}(\alpha)(v) = (\mathbf{D} - \iota_v)(\alpha(v)) \]
and its equivariant curvature 
\[ \mathbf{F}_{\g}(v) = \mathbf{D}_{\g}^2 + \L^{\E}_v = \mathbf{F} +\mu(v)\] 
where $\mathbf{F} = \mathbf{D}^2$ is the curvature of $\mathbf{D}$, and $\mu(v) = \L^{\E}_{v}-[\mathbf{D}, \iota_v]$ is the moment of $v \in \g$ relative to the superconnection $\D$. Then, the \emph{equivariant Chern character of} $\D$ is the equivariant differential form 
\[\ch_{g}(\D)(v) = \mathrm{tr}_s e^{\mathbf{F}_{\g}(v)} = \mathrm{tr}_s e^{\mathbf{F}_{\g}(v)} = \mathrm{tr}_s e^{\mathbf{F}+\mu(v)}\] for every $v \in \g$. \\

Let $\nabla = \nabla^+ \oplus \nabla^-$ be a connection on $\E=\E^+ \oplus \E^-$, with $\nabla^{\pm}$ denoting the parts acting on $E^{\pm}$, and consider its pullback $p^{*}\nabla$ to $p^*E \to TM$ (where $p: TM \to M$ is the canonical projection). Then, consider the superconnection (acting on $\Omega(TM, p^*\E)$):
\[ \mathbf{D} = p^*\nabla + L = p^*\nabla + \sigma_A \otimes 1 + 1 \otimes c(\varphi_x(\xi)). \]

Kasparov's calculations exposed in the previous section suggests the following definition.

\begin{definition} The \emph{equivariant Chern character of the transverse symbol class} $[\sigma_A] \in K_0^G(\mathfrak{S}_\Gamma(M))$ is:
\[ \ch_{\g}[\sigma_A](v) := \dfrac{\ch_{\g}[\D](v)}{\ch_{\g} [W](v)}  \quad ; \quad \forall v \in \g.\]
\end{definition}

\subsection{Example of the circle acting on the complex plane}

We will now do an experiment to check that on the example of the circle action on the complex plane, using the Chern character defined above in the equivariant index formula instead of the Berline--Paradan--Vergne transverse Chern character yields the same result. \\

Let $M=\C$ equipped with the inner product $\langle z_1, z_2 \rangle = \mathrm{Re}(z_1 \ol{z}_2)$ and the usual circle action of $G=S^1$. Let $E=E^{+} \oplus E^{-}$ with $E^{+}=E^{-}=M \times \C$ be the $G$-equivariant (trivial) bundle over $M$ with the action of $G$ on $E$ defined as follows: 
\[ e^{i\theta}\cdot (z, \zeta_+, \zeta_{-}) = (e^{i\theta}z, \zeta_+, e^{i\theta}\zeta_{-}). \] 
Let $A$ be the self-adjoint odd operator with
\[ A: C^\infty(M,E_+) \longrightarrow C^\infty(M,E_-) \quad  ; \quad \left. A\right|_{E_+} = \ol{\partial}_z + z \] 
with symbol
\[ \sigma_A(x,\xi) = 
\begin{bmatrix}
0 & \ol{z} - \i \ol{\xi} \\
z + \i \xi & 0
\end{bmatrix} \] 
Let $z \in M$: for every $\xi \in T_zM$, we have: 
\[ \varphi_z(\xi) = \i z \, \,  \mathrm{ Im}(\ol{z} \xi). \]

The orbital Clifford algebra $\Cl_{\Gamma}(TM)$ is in this case generated, as a module over $C_0(TM)$, by the unit $1 \in \C$ and the section $\varphi_z(\xi)$. So we may embed $\Cliff(\Gamma)$ within the Clifford algebra bundle $W=TM \times \Cliff(\C) \simeq TM \times (\C \oplus \varepsilon\C)$; the first summand being the even part generated by scalars while the second being the odd part generated by $\varepsilon$ such that $\varepsilon^2=1$. The $G$-action on $\Cliff(\Gamma)$ is trivial on the even part (scalars), but the odd part inherits the $G$-action on vectors of $M$ (by rotations). Note that $W$ may also be viewed as a spinor bundle associated to the normal bundle to the fixed point set. \\

According to the decomposition 
\[ p^*E \otimes W = \E^+ \oplus \E^- = (p^*E^+ \otimes W^+ \oplus p^*E^- \otimes W^-) \oplus (p^*E^+ \otimes W^- \oplus p^*E^- \otimes W_+) \]
we have (keeping in mind that tensor products are graded):
\begin{align*}
L = \sigma_A \otimes 1 + 1 \otimes \, c(\varphi_z(\xi)) 
&= \begin{bmatrix}
0 & \ol{z} - \i \ol{\xi} \\
z+\i \xi & 0 
\end{bmatrix} 
\otimes 1 + 
1 \otimes
\begin{bmatrix}
0 & - \i \ol{z} \, \,  \mathrm{ Im}(\ol{z} \xi) \\
\i z \, \,  \mathrm{ Im}(\ol{z} \xi) & 0 
\end{bmatrix} \\
&=\begin{bmatrix}
0 & 0 & -\i \ol{z} \, \,  \mathrm{ Im}(\ol{z} \xi) & \ol{z} - \i \ol{\xi} \\
0 & 0 & z + \i \xi & -\i z \, \,  \mathrm{ Im}(\ol{z} \xi) \\
\i z \, \,  \mathrm{ Im}(\ol{z} \xi) & \ol{z} - \i \ol{\xi} & 0 & 0 \\ z + \i \xi & \i \ol{z} \, \,  \mathrm{ Im}(\ol{z} \xi) & 0 & 0 
\end{bmatrix} 
\end{align*}
A direct determinant calculation shows that $L$ is an elliptic symbol. In the term $1 \otimes c(\varphi_z(\xi))$, use a linear homotopy that transforms the factor $\mathrm{Im}(\ol{z}\xi)$ to $1$, and concatenate this homotopy with the linear one taking $\i z$ to $\i z+\xi$: one checks again by direct calculation that in total, this yields a homotopy of elliptic symbols from $L$ to 
\begin{align*}
\ti{L} & = \begin{bmatrix}
0 & \ol{z} - \i \ol{\xi} \\
z+\i \xi & 0 
\end{bmatrix} 
\otimes 1 + 
1 \otimes
\begin{bmatrix}
0 & - \i \ol{z} + \ol{\xi} \\
\i z + \xi & 0 
\end{bmatrix} \\
& = \begin{bmatrix}
0 & 0 & -\i \ol{z} + \ol{\xi} & \ol{z} - \i \ol{\xi} \\
0 & 0 & z + \i \xi & -(\i z + \xi) \\
\i z + \xi & \ol{z} - \i \ol{\xi} & 0 & 0 \\ z + \i \xi & -(-\i \ol{z} + \ol{\xi}) & 0 & 0 
\end{bmatrix} 
\end{align*}
We now pick the trivial connection $\nabla = d$; using standard transgression formulas and performing the change of variables $u=z+\i \xi$ and $v=iz+\xi$, we can replace the original superconnection on $p^*E \otimes W$ by the following one:
\[ \mathbf{D} = d + \i \ti{L} = \begin{bmatrix}
d & 0 & \i  \ol{v} & \i  \ol{u} \\
0 & d & \i  u & -\i  v \\
\i  v & \i  \ol{u} & d & 0 \\ 
\i  u & -\i  \ol{v} & 0 & d 
\end{bmatrix} \]
with equivariant curvature $\mathbf{F}_{\g} = -(|u|^2+|v|^2)\mathbf{I} + \i \, \mathbf{F}_0 + \mu^{\mathbf{D}}(\i\theta)$, where $\mathbf{F}_0$ and the moment $\mu^{\mathbf{D}}(\i\theta)$ are given by 
\[ \mathbf{F}_0 = \begin{bmatrix}
0 & 0 & d\ol{v} & d\ol{u} \\
0 & 0 & du & -dv \\
dv & d\ol{u} & 0 & 0 \\ 
du & -d\ol{v} & 0 & 0 
\end{bmatrix} \quad ; \quad 
\mu^{\mathbf{D}}(\i\theta) = \begin{bmatrix}
0 & 0 & 0 & 0 \\
0 & 2\i \theta & 0 & 0 \\
0 & 0 & \i \theta  & 0 \\ 
0 & 0 & 0 & \i \theta 
\end{bmatrix} \]
 Using a Duhamel expansion:
\[ e^{\i \mathbf{F}_0+\mu^{\mathbf{D}}(\i\theta)} = e^{\mu^{\mathbf{D}}(\i \theta)} + \sum_{k=1}^4 \int_{\Delta_k} e^{t_1 \mu^{\mathbf{D}}(\i \theta)} \, \i \mathbf{F}_0 \, e^{t_2 \mu^{\mathbf{D}}(\i \theta)} \, \ldots \, \i \mathbf{F}_0 \, e^{t_{k+1} \mu^{\mathbf{D}}(\i \theta)} \, dt_1 \ldots dt_k \]
where $\Delta_k \subset \mathbb{R}^{k+1}$ is the standard $k$-simplex, so one gets in total:
\[ \ch_{\g}(\mathbf{D})(\i\theta) = e^{-(|u|^2+|v|^2)} \left(1+\dfrac{1}{\i \theta}(d\ol{u} du + d\ol{v} dv) - \dfrac{1}{(\i \theta)^2} d\ol{u} dud\ol{v} dv \right)(1-e^{\i\theta})^2. \]
Moreover, an easy calculation shows that
\[ \mathrm{ch}_{\g}(W)(\i\theta) = 1-e^{\i \theta} \]
The Berline--Paradan--Vergne index formula calculates the equivariant index of an elliptic symbol $[\sigma]$ as follows:  
\[ \mathrm{Index}_G [\sigma](e^X) = \int_{TM} \widehat{A}_g(M)^2(X) \, \ti{\mathrm{ch}}_{g}[\sigma](X)  \]
where $\ti{\mathrm{ch}}_{\g}$ is either the Berline--Vergne \cite{BV1, BV2} or the Paradan--Vergne \cite{PV1, PV2, PV3} Chern character, and for $X=\i\theta \in \g$, \[\widehat{A}_{\g}(M)^2(\i\theta) = \det\left(\frac{R^{TM}_{\g}/2}{\sinh(R^{TM}_{\g}/2)}\right) = \frac{(\i \theta)^2 \, e^{\i \theta}}{(1-e^{\i \theta})^2}. \] with $R^{TM}_{\g}$ is the equivariant curvature of a connection on $p^*TM$. \\

In view of the above index formula, we evaluate the following quantity, with is the same as the right-hand-side with our Chern character being used instead:
\[ \int_{TM} \widehat{A}_{\g}(M)^2(\i \theta) \, \ch_{\g}[\sigma_A](\i \theta) \]
In addition, the integral over $TM$ selects only the term of $\ch_{\g}[\sigma_A]$ attached to the volume form $d\ol{u} dud\ol{v} dv$. Overall, we get
\[ \int_{TM} \widehat{A}_{\g}(M)^2(\i \theta) \, \mathrm{ch}_{\g}[\sigma_A](\i \theta) = - \dfrac{e^{\i \theta}}{1 - e^{\i \theta}} = -\sum_{n\geq 1} e^{\i n \theta}. \]
which is consistent with the calculations of Berline--Paradan--Vergne, and also with the result of Kasparov outlined in the previous section. The result should be viewed in the sense of distributions on $\g$. 

\subsection{Example of the zero operator}

Consider the natural circle action of $G=S^1$ onto $M=S^1$, and let $\mathbf{0}: C^{\infty}(S^1,\C) \to C^{\infty}(S^1,\{0\}) = \{0\}$ be the zero-operator with symbol $\sigma(x,\xi)=0 \in C^{\infty}(TS^1, \Hom(S^1 \times \C, S^1 \times \{0\})$. Using the parametrization $TS^1 = S^1 \times \mathbb{R} = \{(e^{\i\theta}, \xi) \, ; \, \theta \in [0, 2\pi]\}$, we have 
\[\varphi_\theta(\xi) = \xi \in T_\theta S^1 \]
that we may identify with the covector $\xi d\theta \in T_{\theta}^{*} S^1$, which is basically the Liouville canonical 1-form (when lifting it to $T^*(T^*S^1)$). The fundamental vector field is 
\[ \zeta(X) = X \partial_{\theta} \quad ; \quad \forall X \in \g. \]
The equivariant superconnection associated to the symbol of the zero-operator is 
\[ \D_{\g} = d - \iota_{\zeta} + \i \omega \]
which has equivariant curvature
\[ \mathbf{F}_g(X) = \D_{\g}^2 = \i(d\theta \, d\xi - X\xi) \quad ; \quad \forall X \in \g. \] 
In this context, $\widehat{A}_{\g}(S^1) = 1$ and 
\[ \ch_{\g}(\sigma) = e^{\i d\theta d\xi} e^{\i \xi X} = (1+\i d\theta d\xi)e^{\i\xi X}. \]
So the index of the zero operator is given by
\[ \mathrm{Index}(\sigma)(e^{X}) = \dfrac{1}{2\pi\i}\int_{S^1 \times \mathbb{R}} e^{\i\xi X} \, \i d\theta d\xi = \int_{\mathbb{R}} e^{\i \xi X} \, d\xi = \delta(0). \]
where $\delta$ is the Dirac delta function at $0 \in \g$.  \\

%\begin{enumerate}
%\item[\bullet] $\rho_z(\i \theta) = \dfrac{d}{dt}\right|_{t=0} e^{\i t \theta} = \i \theta z$. 
%\item[\bullet]$\rho_z^{*}(\xi) = \i z \, \,  \mathrm{ Im}(\ol{z} \xi)$. 
%\end{enumerate}

%\begin{remark}
%It is also possible to do this construction without reduction to a torus action, but one should be cautious about the fact that $p^*E \otimes \Cl_{\Gamma}(TM)$ is not a vector bundle (e.g there are jumps in its `rank'). Nonetheless, as $\Cl_{\Gamma}(M)$ is a space of sections of the action Lie algebroid, generalizing the Chern character to this context is possible.   
%\end{remark}

%%%%%%%%%%%%%%%%%%%%%%%%%%%%%%%%%%%%%%%%
%%%%%%%%%%%% References %%%%%%%%%%%%%%%

\providecommand{\bysame}{\leavevmode\hbox to3em{\hrulefill}\thinspace}
\providecommand{\MR}{\relax\ifhmode\unskip\space\fi MR }
% \MRhref is called by the amsart/book/proc definition of \MR.
\providecommand{\MRhref}[2]{%
  \href{http://www.ams.org/mathscinet-getitem?mr=#1}{#2}
}
\providecommand{\href}[2]{#2}

\end{document}